\documentclass[11pt,reqno]{amsart}

\usepackage{amssymb, amsmath,euscript}

\newtheorem{theorem}{Theorem}
\newtheorem{lemma}[theorem]{Lemma}
\newtheorem{proposition}[theorem]{Proposition}
\newtheorem{corollary}[theorem]{Corollary}
\newtheorem{example}[theorem]{Example}
\newtheorem{definition}[theorem]{Definition}
\newtheorem{remark}[theorem]{Remark}

\renewcommand{\phi}{\varphi}

\DeclareMathSymbol{\subsetneqq}{\mathrel}{AMSb}{"24}
\DeclareMathSymbol{\supsetneqq}{\mathrel}{AMSb}{"25}
\DeclareMathSymbol{\varnothing}{\mathord}{AMSb}{"3F}

\DeclareMathOperator{\supp}{supp}
\DeclareMathOperator{\diam}{diam}
\DeclareMathOperator{\cov}{Cov}
\def\ee{\varepsilon}
\def\e{\mathrm{e}}

\def\T{\mathbb{T}}
\def\Z{\mathbb{Z}}

\def\cA{\EuScript{A}}

\def\uR{\underline{R}}
\def\oR{\overline{R}}
\def\udmu{\underline{d}_\mu}
\def\odmu{\overline{d}_\mu}
\renewcommand{\phi}{\varphi}

\def\RR{\mathbb{R}}

\def\lb{[}
\def\rb{]}

\begin{document}
\date{November 2004}
\title{Recurrence rate in rapidly mixing dynamical systems}
\author{Benoit Saussol}
\email{benoit.saussol@u-picardie.fr}
\urladdr{http://lamfa.u-picardie.fr/saussol}
\address{LAMFA - CNRS UMR 6140\\ Universit\'e de Picardie Jules Verne\\ 33 rue St Leu, 80039 Amiens cedex 1\\ France}

\begin{abstract}
For measure preserving dynamical systems on metric spaces we study the time needed by a typical orbit to return back close to its starting point.
We prove that when the decay of correlation is super-polynomial the recurrence rates and the pointwise dimensions are equal.
This gives a broad class of systems for which the recurrence rate equals the Hausdorff dimension of the invariant measure.
\end{abstract}

\maketitle

\section{Introduction}

\subsection{Decay of correlations}

Let $(X,f,\mu)$ be a measure preserving dynamical system. Recall that the system is said to be \emph{mixing}
if for any functions $\phi$, $\psi$ in $L^2$ the covariance 
\begin{equation}\label{eq:doc}
\cov(\phi\circ f^n,\psi) := \int \phi\circ f^n\psi d\mu - \int \phi d\mu \int \psi d\mu\to 0
\quad\text{as $n\to\infty$.}
\end{equation}
The decay of the correlation function is, in great generality, arbitrarily slow. 
The notion of rapid mixing needs a little more structure. 

Assume that $X$ is a metric space with metric $d$, and consider the space $\text{Lip}(X)$ of real Lipschitz functions on $X$.
For many dynamical systems an upper bound for~\eqref{eq:doc} of the form
$\|\phi\| \|\psi\| \theta_n$ has been computed, where $\theta_n\to 0$ with some rate, and $\|\cdot\|$ is a norm on a space of functions with some regularity. Without loss of generality we are considering in this paper the rate of decay of correlations for Lipschitz observables\footnote{For example an immediate approximation argument allows  easily to go from Holder or class $C^k$ to Lipschitz.}.

A broad class of systems enjoy exponential decay of correlations. The main result of the paper (Theorem~\ref{thm:B}) applies to systems with super-polynomial decay of correlation. 
This includes for example Axiom~A systems with equilibrium states, hyperbolic systems with singularities with their SBR measures such as those considered by chernov in~\cite{ch}, many systems with a Young tower~\cite{yo1,yo2}, expanding maps with singularities such as in~\cite{sa}, some non-uniformly expanding maps~\cite{al}, etc.
The main reference for these questions is certainly the book by Baladi~\cite{ba}. The reader will also find in the review by Luzzatto~\cite{lu} an exposition of the recent methods for non-uniformly expanding systems and an extensive bibliography on this active field.

\subsection{Recurrence rate and dimensions}

The return time of a point $x\in X$ under the map $f$ in its $r$-neighborhood is
\begin{equation*}
\tau_r(x)=\inf\{n\ge1\colon d(f^nx,x)<r\}.
\end{equation*}
We are interested in the behavior as $r\to0$ of the return time.
We define the recurrence rate as the limits
\begin{equation*}
\uR(x)=\liminf_{r\to0}\frac{\log\tau_r(x)}{\log(1/r)}
\quad\text{and}\quad
\oR(x)=\limsup_{r\to0}\frac{\log\tau_r(x)}{\log(1/r)}.
\end{equation*}
Whenever $\uR(x)=\oR(x)$ we denote by $R(x)$ the value of the limit.

{} From now on we assume that $X$ is a finite dimensional Euclidean space.
Denote by $HD(Y)$ the Hausdorff dimension of a set $Y\subset X$.
We define the Hausdorff dimension of a probability measure $\mu$ by
\begin{equation*}
HD(\mu) = \inf \{HD(Y)\colon \mu(Y)=1\}
\end{equation*}
We also define a local version of the dimension, namely
\begin{equation}
\udmu(x)=\liminf_{r\to0}\frac{\log\mu(B(x,r))}{\log r}
\quad\text{and}\quad
\odmu(x)=\limsup_{r\to0}\frac{\log\mu(B(x,r))}{\log r}
\end{equation}
It is well known that the Hausdorff dimension satisfies the relation 
\begin{equation}\label{eq:dim}
HD(\mu)  = \text{ess-sup } \udmu.
\end{equation}

Barreira and Saussol established in~\cite{bs1} the following relation 
\begin{proposition}\label{pro:bs}
Let $f$ be a measurable map and $\mu$ be an invariant measure for $f$.
The recurrence rates are bounded from above by the pointwise dimensions~:
\[
\uR\le \udmu
\quad\text{and}\quad
\oR\le \odmu
\quad\mu\text{-a.e.}
\]
\end{proposition}
We refer to the works by Boshernitzan~\cite{bo} and Ornstein and Weiss~\cite{ow} for pioneering related results.

In this paper we are giving conditions under which the opposite inequalities will hold, establishing the equalities
\begin{equation}\label{eq:cequonveux}
\uR=\udmu
\quad\text{and}\quad
\oR= \odmu
\quad\mu\text{-a.e.}
\end{equation}

\subsection{Statement of the results}

\begin{definition} We say that $(X,f,\mu)$ has super-polynomial decay of correlations  
if we have
\begin{equation}\label{eq:ddc}
\left|\int \phi\circ f^n \psi d\mu - \int \phi d\mu\int \psi d\mu\right|
\le \|\phi\| \|\psi\| \theta_n
\end{equation}
with $\lim_n\theta_n/n^p=0$ for all $p>0$, where $\|\cdot\|$ is the Lipschitz norm.

We say that the \emph{local} decay of correlations is super-polynomial
if there exists a partition (modulo $\mu$) into open sets $V_i$ and sequences $\theta^i_n$ such that \eqref{eq:ddc} holds whenever $\supp \phi\subset V_i$ and $\supp \psi\subset V_i$, where $\lim_n\theta_n^i/n^p=0$ for all $p>0$.
\end{definition}

The main result of the paper is the following.
\begin{theorem}\label{thm:B}
Let $(X,f,\mu)$ be a measure preserving dynamical system.
If the entropy $h_\mu(f)>0$, $f$ is Lipschitz (or piecewise Lipschitz with finite average Lipschitz exponent ; see Definition~\ref{def:lip}) and the (local) decay of correlation is super-polynomial then
\[
\uR=\udmu
\quad\text{and}\quad
\oR=\odmu
\quad\mu\text{-a.e.}
\]
\end{theorem}
We postpone the proof at the end of Section~\ref{sec:nonzero}.
This extends some results by Barreira and Saussol in~\cite{bs1,bs2}, including the case of Axiom~A systems with equilibrium states. 
The theorem also applies to loosely Markov dynamical systems and we recover Urbanski's result in~\cite{ur}.
The hypotheses in Theorem~\ref{thm:B} are satisfied in a number of systems such as those already quoted in the introduction. 
All these systems have in common some hyperbolic behavior. We now give an example of a relatively different nature, due to the possibility of zero Lyapunov exponents, where one can still apply Theorem~\ref{thm:B}.

\begin{example}[Ergodic toral automorphisms] 
Recall that any matrix $A\in \text{Sl}(k,\Z)$ (i.e. the entries of $A$ are in $\Z$ and $|\det A|=1$) gives rise to an automorphism $f$ of the torus $\T^k$ by $f(x)=Ax \mod \Z^k$ which preserves the Lebesgue measure. The map $f$ is ergodic if and only if the matrix $A$ has no eigenvalue root of unity.
Lind's established~\cite{li} the exponential decay of correlations (using the algebraic nature and Fourier transform) which is more than enough to apply Theorem~\ref{thm:B} and get
\[
R(x)=k\quad\text{for Lebesgue a.e. $x\in\T^k$}.
\]
for any ergodic automorphism of the torus, even \emph{non-hyperbolic}. 
\end{example}

Let $f$ be a diffeomorphism of a compact manifold $M$ and $\mu$ be an invariant measure. 
By Oseledec's multiplicative ergodic Theorem the Lyapunov exponents
\[
\lambda(x,v)=\lim\frac1n \log |d_xf^nv|
\]
are well defined for all nonzero $v \in T_xM$ for a.e. $x\in M$.
Recall that a measure $\mu$ is said to be hyperbolic if none of its Lyapunov exponents are zero. Barreira, Pesin and Schmeling~\cite{bps} prove the following.
\begin{proposition}\label{pro:1}
Let $f$ be a diffeomorphism of a compact manifold  and $\mu$ be an ergodic hyperbolic measure. Then we have
\[
\udmu=\odmu=HD(\mu)
\quad\mu\text{-a.e.}
\]
\end{proposition}

The case of an hyperbolic measure with zero entropy is completely understood.
\begin{proposition}\label{pro:0}
let $f$ be a diffeomorphism of a compact manifold and $\mu$ be an hyperbolic invariant measure. 
If $h_\mu(f)=0$ then $R=0=HD(\mu)$ $\mu$-a.e.
\end{proposition}

\begin{proof}
Barreira and Saussol established in~\cite{bs1} the inequality $\oR\le \odmu$ $\mu$-a.e.
and it follows from Ledrappier and Young's work~\cite{ly} that $HD(\mu)=0$ if $h_\mu(f)=0$, which allows to conclude by Proposition~\ref{pro:1}.
\end{proof}

\begin{corollary}
Let $f$ be a diffeomorphism of a compact manifold and $\mu$ be an hyperbolic measure with super-polynomial rate of decay of correlation. Then we have
\[
R=HD(\mu)
\quad\mu\text{-a.e.}
\]
\end{corollary}
\begin{proof}
If the entropy is zero then this is the content of Proposition~\ref{pro:0}. 
If the entropy is non-zero then this is the content of Theorem~\ref{thm:B}.
\end{proof}
We point out that in the case of interval maps with nonzero Lyapunov exponent, Saussol, Troubetzkoy and Vaienti prove that $R=HD(\mu)$ $\mu$-a.e. for ergodic measures, under very weak regularity conditions~\cite{stv}. See  Remark~\ref{rem:confbound}-(i) for related results.
\\

We now give a sketch of the strategy adopted in this paper.

Theorem~\ref{thm:A} states that under sufficiently rapid mixing the recurrence rates equal the pointwise dimensions a.e. on the set where $\uR>0$.
Indeed, mixing implies that
$\mu(B\cap f^{-n}B)\to \mu(B)^2$ as $n\to\infty$. Thus we have $\mu(B\cap f^{-n}B)\le 2\mu(B)^2$ for large $n$.
If now we consider the set $B\cap f^{-n}B\cap f^{-n-1}B\cap\cdots\cap f^{-n-\ell}B$ then its measure is bounded by $2\ell\mu(B)^2$. If $\ell\le\mu(B)^{-1+\ee}$ then we get that the proportion of points inside $B$ that never enter in $B$ in the time interval $[n,n+\ell]$ is bounded by $2\mu(B)^\ee$. Using the decay of correlations we are able to prove that this last statement is true for $n$ of the order $\diam(B)^{-\delta}$ for some small $\delta>0$, whenever $B$ is a ball.
This is what we call the \emph{long fly property}. A Borel Cantelli argument then shows that typical points do have long flies (see Lemma~\ref{lem:longfly} for precise statement). If in addition we also have $\uR>\delta$ then it immediately shows that the return time into small neighborhoods $B$ cannot be much less (at an exponential scale) than $\mu(B)^{-1}$, establishing Equation~\eqref{eq:cequonveux}.

On the other hand, for systems which are not too wild (e.g. finite Lyapunov exponents, see Lemma~\ref{lem:reasen}) and with nonzero metric entropy, a symbolic coding (see Lemma~\ref{lem:largint}) allows to use Ornstein-Weiss' theorem on repetition time of symbolic sequences to prove that the return time of a typical point in a ball $B$ is not less than $\diam(B)^{-\delta}$ ; see Lemma~\ref{lem:R>0}.

The structure of the paper is as follows. 
We state and prove in Section~\ref{sec:rapid} the core result, Theorem~\ref{thm:A}. 
In Section~\ref{sec:nonzero} we provide some conditions under which the recurrence rate is nonzero.

\section{Rapid mixing implies long flies}\label{sec:rapid}

\begin{theorem}\label{thm:A}
Assume that the local rate of decay of correlations is super-polynomial. Then on the set $\{\uR>0\}$ we have
\[
\uR=\udmu \quad\text{and}\quad \oR=\odmu
\quad\mu\text{-a.e.}
\]
\end{theorem}

\begin{proof}
By Proposition~\ref{pro:bs} we know that $\uR\le\udmu$ and $\oR\le\odmu$.
Furthermore, the first inequality implies that $\{\uR>a\}\subset\{\udmu>a\}$ $\mu$-a.e.
But on the set $\{\uR>a\}$ we have $\tau_r(x)\ge r^{-a}$ provided $r$ is sufficiently small. 
By Lemma~\ref{lem:longfly} below with $\delta=a$ and $\ee>0$ we get that $\tau_r(x)\ge \mu(B(x,r))^{-1+\ee}$ provided $r$ is sufficiently small, for $\mu$-a.e. $x\in \{\uR>a\}$.  Thus 
$\uR\ge (1-\ee)\udmu$ and $\oR\ge(1-\ee)\odmu$ $\mu$-a.e. on $\{\uR>a\}$.
The conclusion follows by taking $\ee>0$ arbitrary small.
\end{proof}

The following lemma expresses that the orbit of a typical point has the long fly property. 
\begin{lemma}\label{lem:longfly}
Let $X_a=\{\udmu>a\}$ for some $a>0$. For any $\delta,\ee>0$, for $\mu$-a.e. $x\in X_a$ there exists $r(x)>0$ such that 
for any $r\in(0,r(x))$ and any integer $n$ in $[r^{-\delta},\mu(B(x,r))^{-1+\ee}]$ we have $d(f^nx,x)\ge r$.
\end{lemma}
\begin{proof}
For clarity we assume that the (global) rate of decay of correlation is super-polynomial. The obvious modifications in the proof would consits essentially in considering separately each sets $G\cap \{x\in V_i\colon d(x,\partial V_i)>\nu\}$ for arbitrarily small $\nu>0$ in place of the unique set $G$ defined below.
 
Let $D=\dim(X)$.
Fix $b>0$, $c=a\ee/3$ and consider for $r_0>0$ the set $G = G_1 \cap G_2\cap G_3$ where
\[
\begin{split}
G_1 &=\{x\in X_a\colon\forall r\le r_0,  \mu(B(x,r))\le r^a\} \\
G_2 &=\{x\in X\colon\forall r\le r_0, \mu(B(x,r))\ge r^{D+b}\} \\
G_3 &=\{x\in X\colon\forall r\le r_0,\mu(B(x,r/2))\ge \mu(B(x,4r))r^{c}\}.
\end{split}
\]
We claim that $\mu(G)\to\mu(X_a)$ as $r_0\to0$. 
Indeed, by definition of the lower pointwise dimension we have $\mu(G_1)\to\mu(X_a)$.
In addition since $\odmu\le D$ a.e. we have $\mu(G_2)\to1$ and since $X$ is Euclidean the measure $\mu$ is weakly diametrically regular (see Lemma~1 in \cite{bs1}), thus $\mu(G_3)\to1$ as well.
Let $r\le r_0$ and define the set
\[
A_\ee(r)=\{y\in X\colon \exists n\in \lb r^{-\delta},\mu(B(y,3r))^{-1+\ee}\rb, d(f^ny,y)<r\}.
\]
Let $x\in G$. By the triangle inequality we get the inclusions
\[
\begin{split}
B(x,r)\cap A_\ee(r)
&\subset \{y\in B(x,r)\colon \exists n\in \lb r^{-\delta},\mu(B(x,2r))^{-1+\ee}\rb, d(f^ny,x)<2r\}\\
&\subset \bigcup_{r^{-\delta}\le n \le \mu(B(x,2r))^{-1+\ee}} B(x,r)\cap f^{-n}B(x,2r). 
\end{split}
\]
Let  $\eta_r\colon[0,\infty)\to\RR $ be the $r^{-1}$-Lipschitz map such that $1_{[0,r]}\le\eta_r\le1_{[0,2r]}$ and set $\phi_{x,r}(y)=\eta_r(d(x,y))$. Clearly $\phi_{x,r}$ is also $r^{-1}$-Lipschitz.
By the assumption on the decay of correlation we obtain
\[
\begin{split}
\mu(B(x,r)\cap f^{-n} B(x,2r)) 
&\le
\int \phi_{x,2r}\phi_{x,2r}\circ f^n d\mu \\
&\le
 \|\phi_{x,2r}\|^2 \theta_n +\left(\int\phi_{x,2r}d\mu\right)^2\\
&\le 
 r^{-2} \theta_n +\mu(B(x,4r))^2.
\end{split}
\]
Choose $p>1$ such that $\delta(p-1)-2\ge D+2b$ and take $r_0$ so small that $n\ge r_0^{-\delta}$ implies $\theta_n\le (p-1) n^{-p}$.
Since $\sum_{n\ge q} n^{-p}\le \frac{1}{p-1} q^{1-p}$ we obtain 
\[
\begin{split}
\mu( B(x,r)\cap A_\ee(r) ) 
&\le r^{\delta(p-1)-2} + \mu(B(x,2r))^{-1+\ee}\mu(B(x,4r))^2\\
&\le
\mu(B(x,r/2)) \left(r^b+ r^{\ee a-2c}\right).
\end{split}
\]
Let $B\subset G$ be a maximal $r$-separated set\footnote{that is if $x\neq x'\in B$ then $d(x,x')\ge r$ and maximal in the sense that for any $y\in G$ there exists $x\in B$ such that $d(x,y)<r$.}. Since $(B(x,r))_{x\in B}$ covers $G$ we have 
\[
\begin{split}
\mu(G\cap A_\ee(r)) 
&\le \sum_{x\in B} \mu( B(x,r) \cap A_\ee(r) )\\
&\le \sum_{x\in B} \mu(B(x,r/2)) (r^b+ r^{\ee a-2c})\\
&\le r^b+ r^{\ee a-2c}
\end{split}
\]
since by the balls $(B(x,r/2))_{x\in B}$ are disjoints.
This implies that 
\[
\sum_m \mu(A_\ee(\e^{-m}))<\infty,
\]
thus by Borel-Cantelli Lemma we obtain that for $\mu$-a.e. $y\in G$ there exists $m(y)$ such for every $m>m(y)$ 
there exists no $n\in\lb \e^{-\delta m},\mu(B(y,3\e^{-m}))^{-1+\ee}\rb$ such that $d(f^ny,y)<\e^{-m}$. 
By weak diametric regularity (and changing slightly if necessary the values of $\ee$ and $\delta$), this proves the lemma.
\end{proof}

\begin{remark}\label{rem:ddcopt}
Observe that we only use that the decay of correlation is at least $n^{-p}$ for some $p>\frac{D+2}{\delta}+1$.
If in addition \eqref{eq:ddc} holds with the first norm $\|\phi\|$ taken to be the $L^1(\mu)$ norm (e.g. expanding maps) then $p>\frac{D+1}{\delta}+1$ suffices.
\end{remark}

\section{Non-zero recurrence rate}\label{sec:nonzero}

We proceed now to find conditions under which the recurrence rate does not vanish.
Denote by $\xi(x)$ the unique element of a partition $\xi$ containing the point $x$ and 
by $\xi^n=\xi\vee f^{-1}\xi\vee\cdots\vee f^{-n+1}\xi$ the dynamical partition, for any integer $n$ .

\subsection{Coding by symbolic systems : partitions with large interior}

\begin{definition}
We say that a partition $\xi$ has large interior if for $\mu$-a.e. $x$ there exists $\chi=\chi(x)<\infty$ such 
that $B(x,\e^{-\chi n})\subset\xi^n(x)$ for all $n$ sufficiently large.
\end{definition}

Next lemma, which proof is fairly simple, is the key-observation which gives to Theorem~\ref{thm:A} all its interest.

\begin{lemma}\label{lem:R>0}
If there exists a partition with large interior and nonzero entropy then $\uR>0$ $\mu$-a.e.
\end{lemma}
\begin{proof}
Let $\xi$ be such a partition. Define 
\[
R_n(x,\xi)=\min\{k>0\colon f^kx\in\xi^n(x)\}.
\]
Ornstein and Weiss \cite{ow} prove that if $\xi$ is a finite partition with entropy $h_\mu(f,\xi)$ then
\[
\lim_{n\to\infty} \frac1n \log R_n(x,\xi)=h_\mu(f,\xi)
\quad\mu\text{-a.e.}
\]
Since $\xi$ has large interior, for $\mu$-a.e. $x\in X$ there exists a number $\chi=\chi(x)$ such that $B(x,e^{-\chi n})\subset \xi^n(x)$. 
Thus 
\[
\uR(x) = \liminf_{n\to\infty} \frac{\log \tau_{\e^{-n}}(x)}{n\chi(x)}
\ge \liminf_{n\to\infty} \frac{\log R_n(x,\xi)}{n\chi(x)} = \frac{h_\mu(f,\xi)}{\chi(x)}>0
\quad\mu\text{-a.e.}
\]
\end{proof}

Combining Lemma~\ref{lem:R>0} and Theorem~\ref{thm:A} we get that if we have local super-polynomial decay of correlations and a partition of positive entropy with large interior then $\uR=\udmu$ and $\oR=\odmu$. 
The rest of the section consists in finding sufficient conditions for the existence of such a partition.

\subsection{Reasonable dependence on initial condition}

\begin{definition}\label{def:reasen}
We say that a system $(X,f,\mu)$ is reasonably sensitive if for $\mu$-a.e. $x$ there exists $\gamma,\lambda>0$ such that $f^n$ is $\e^{\lambda n}$-Lipschitz on the ball $B(x,\e^{-\gamma n})$ for all $n$ sufficiently large.
\end{definition}

\begin{lemma}\label{lem:largint}
If the system $(X,f,\mu)$ is reasonably sensitive and the entropy $h_\mu(f)>0$ then there exists a partition with large interior and nonzero entropy.
\end{lemma}
\begin{proof}
{\bf Claim :} For any $x\in X$, $s>0$ there exists $\rho\in(s,2s)$ such that 
\begin{equation}\label{eq:cl1}
\mu( \{y\in X\colon \rho-4^{-n}s<d(x,y)<\rho+4^{-n}s\} ) \le \frac{1}{2^{n-1}} \mu( B(x,2s) ).
\end{equation}
Indeed, let $m$ be the measure on the interval $(0,2)$ defined by $m([0,t))=\mu(B(x,st))$.
We construct a sequence of open intervals $I_n$ starting from $I_0=(1,2)$.  If $I_n$ is an interval of length $4^{-n}$ we divide it into 4 pieces of equal length and choose $I_{n+1}$ the left of the right central piece of smallest measure. We have $m(I_{n+1})\le\frac12 m(I_n)$. $I_n$ is a decreasing sequence of intervals with $\overline I_{n+1}\subset I_n$ thus $\cap_n I_n$ contains one point, say $\bar\rho$. Since $\bar\rho\in I_n$ we have $\bar\rho\pm 4^{-n}\in I_{n-1}$ thus 
$m((\bar\rho-4^{-n},\bar\rho+4^{-n}))\le m(I_{n-1})\le \frac1{2^{n-1}} m(I_0)$.
Proving the claim with $\rho=s\bar\rho$.
\\

Fix $s>0$ so small that any partition made by sets of diameter less than $2s$ has nonzero entropy.
Choose a maximal $s$-separated set $E$. For any $x\in E$ take $\rho_x\in (s,2s)$ such that \eqref{eq:cl1} in the claim holds.
Let $E=\{x_1,x_2,\ldots\}$ be an enumeration of the (at most) countable set $E$. Put $B_i=B(x_i,\rho_{x_i})$ and define
$Q_1=B_1$, $Q_2=B_2\setminus Q_1$, $Q_3=B_3\setminus (Q_1\cup Q_2)$, $\ldots$
By maximality the collection of sets $\xi=\{Q_1,Q_2,\ldots\}$ is a partition of $X$ (modulo $\mu$) and since $\partial\xi\subset \cup_i\partial B_i$ we get
\[
\begin{split}
\mu(\{x\in X\colon d(x,\partial \xi) <4^{-n}s\})
&\le \mu(\cup_i \{x\in X\colon \rho_{x_i}-4^{-n}<d(x_i,x)<\rho_{x_i}+4^{-n}\})\\
&\le \frac{1}{2^{n-1}} \sum_i \mu(B(x_i,2s)).
\end{split}
\]
Since the $x_i$ are $s$-separated and $X$ is Euclidean there are at most $c(X)=c(\dim X)$ balls of radius $2s$ that can intersect, thus the last sum is bounded by $\frac{c(X)}{2^{n-1}}$. 
This proves that for some constants $a,c>0$ and all $\ee>0$
\[
\mu(x\in X\colon d(x,\partial\xi)<\ee)<c\ee^a.
\]
Thus for any $b>0$ we have by the invariance of $\mu$ 
\[
\sum_n \mu(\{x\in X\colon d(f^nx,\xi)<\e^{-bn}\}) \le \sum_nc\e^{-abn}<\infty.
\]
This implies by Borel-Cantelli Lemma that for $\mu$-a.e. $x$ there exists $n(x)<\infty$ such that $d(f^nx,\partial\xi)\ge\e^{-bn}$, hence $B(f^nx,\e^{-bn})\subset \xi(f^nx)$, for any $n\ge n(x)$.
Taking $c(x)\in(0,1)$ sufficiently small we have $B(f^nx,c(x)\e^{-bn})\subset \xi(f^nx)$ for all integer $n$.

Fix $x\in X$ where the reasonable sensitovoty condition holds. Without loss of generality, and changing if necessary $c(x)$ into a smaller constant we assume that $f^n$ is $\e^{\lambda n}$-Lipschitz on the ball $B(x,c(x)\e^{-\gamma}n)$ for all integer $n$ and that $\lambda>\gamma+b$.

We show then by induction that $B(x,c(x)\e^{-\lambda n})\subset \xi^k(x)$ for any $k\le n$.
Indeed, this is trivially true for $k=1$, and if this holds for some $k\le n-1$ then we have 
\[
f^k(B(x,c(x)^2\e^{-\gamma n})\subset B(f^kx,c(x)\e^{\lambda k-\gamma n})\subset B(f^kx,\e^{-bn})\subset\xi(f^kx).
\] 
Hence 
$B(x,c(x)^2\e^{-\gamma n})\subset \xi^{k+1}(x)$.
\end{proof}

We finally provide a sufficient condition for reasonable sensitivity.

\begin{definition}\label{def:lip}
If there exists a partition $\cA$ (modulo $\mu$) into open sets such that on each $A\in \cA$ the map $f$ is Lipschitz with constant $L_f(A)$ and the singularity set $\partial \cA=\cup_{A\in\cA}\partial A$ is such that
$\mu( \{x\in X\colon d(x,\partial\cA)<\epsilon\}) \le c\epsilon^a$ for some constants $c>0$ and $a>0$ then we say that $f$ is piecewise Lipschitz with average Lipschitz exponent $\log L_f=\int \log^+ L_f(\cA(x))d\mu(x)=\sum_{A\in\cA}\log^+ L_f(A)\mu(A)$.
\end{definition}

\begin{lemma}\label{lem:reasen}
If $f$ is Lipschitz, or piecewise Lipschitz with finite Lispchitz exponent
then $(X,f,\mu)$ is reasonably sensitive.
\end{lemma}
\begin{proof}
We prove the piecewise case, the other one is obvious.
Let $\lambda>\log L_f$.
By the Birkhoff Ergodic Theorem, for $\mu$-a.e. $x$ there exists $m(x)$ such that 
\[
L_f(\cA(x))L_f(\cA(fx))\cdots L_f(\cA(f^{n-1}x)) \le \e^{\lambda n}
\] 
for all $n\ge m(x)$. Replacing if necessary the upper bound by $\e^{\lambda n}/c(x)$ for some constant $c(x)\ge1$ the inequality will hold for any integer $n$.
Proceeding as in the last part of the proof of Lemma~\ref{lem:largint} we get that 
for any $b>0$, changing $c(x)$ if necessary, we have $B(f^nx,c(x)\e^{-bn})\subset \cA(f^nx)$ for any integer $n$.
We then conclude similarly that $B(x,c(x)^2\e^{-bn}\e^{-\lambda n})\subset \cA^n(x)$.
This concludes the proof taking $\gamma=b+\lambda$.
\end{proof}

The proof of Theorem~\ref{thm:B} follows now easily from the preceding results.
\begin{proof}[Proof of Theorem~\ref{thm:B}]
By Lemma~\ref{lem:reasen} the map is reasonably sensitive. 
This implies by Lemma~\ref{lem:largint} the existence of a partition with large interior.
By Lemma~\ref{lem:R>0} we find that $\uR>0$ a.e. and the conclusion follows from Theorem~\ref{thm:A}.
\end{proof}

\begin{remark}\label{rem:confbound}
(i)
We remark that if $f$ is $C^1$ on a compact manifold, or more generally if $f$ is piecewise $C^{1+\alpha}$ with reasonable singularity set such as in \cite{ks}, 
then the exponents $\lambda$ and $\gamma$ in Definition~\ref{def:reasen} can be taken arbitrarily close to the largest Lyapunov exponent\footnote{to see this, consider a Lyapunov chart whose local chart at $x$ has a diameter $\rho(x)$, where $\rho$ is $\eta$-slowly varying. A choice like $\lambda=\lambda_\mu^+ +2\eta$ and $\gamma=\lambda+\eta$ would do the job.} $\lambda_\mu^+$.
Thus the exponent $\chi$ in Lemma~\ref{lem:R>0} may also be taken arbitrarily close to $\lambda_\mu^+$. This readily implies that $\uR\ge h_\mu/\lambda_\mu^+$. This is optimal in dimension one or more generally for conformal maps, where under mild assumptions we have $HD(\mu)=h_\mu/\lambda_\mu$.
\\

(ii)
Combining the above observation with Remark~\ref{rem:ddcopt} shows that the assumption on the super-polynomial decay of correlations in Theorem~\ref{thm:A} may be reduced to a decay at a rate $n^{-p}$ for some $p>\frac{D+2}{h_\mu}\lambda_\mu^+ +1$.
\end{remark}

\end{document}